\documentclass[12pt]{article}
\usepackage{mathrsfs}
\usepackage{epic,eepic,epsf,epsfig}
\usepackage{amsfonts,srcltx,mathrsfs}
\usepackage{tikz}
\usetikzlibrary{calc}
\usepackage{amsmath}
\usepackage{amssymb}
\usepackage{amsbsy}
\usepackage{float}
\usepackage{graphicx}
\usepackage{amsfonts}
\usepackage{color}

\usepackage{url}
\newtheorem{lem}{Lemma}[section]
\newtheorem{theorem}[lem]{Theorem}

\newtheorem{case}{Case}
\newtheorem{prop}[lem]{Proposition}

 \def\b{\beta}   
\def\r{\rho} \def\s{\sigma}

\def\PGL{\hbox{\rm PGL}}
\def\PGaL{\hbox{\rm P\mbox{$\Gamma$}L}}

 \def\lg{\langle} \def\rg{\rangle}
\def\nd{\mathrel{\bigm|\kern-.7em/}}

\def\f{\noindent}

\def\PSL{\hbox{\rm PSL}}
\def\PGL{\hbox{\rm PGL}}

 \def\PGL{\hbox{\rm PGL}}
\def\Sz{\hbox{\rm Sz}}

\def\P{\mathcal{P}}

\def\Aut{\hbox{\rm Aut}}

\def\Aut{\hbox{\rm Aut}}

\def\PGL{\hbox{\rm PGL}}

\def\demo{\f {\bf Proof.}\hskip10pt}

\newcommand{\qed}{\mbox{\raisebox{0.7ex}{\fbox{}}} \vspace{4truemm}}

\begin{document}

\title{A family of regular polytopes of order $4p^m$ with type $\{p, 2p\}$}

\author{ \\ Dong-Dong Hou, Ting-Ting Kong, Hai-Peng Qu\\
{\small Department of Mathematics, Shanxi Normal University}\\
{\small Linfen, Shanxi
041004, P.R. China}\\
}

\date{}
\maketitle

\footnotetext{Author.E-mails: holderhandsome$@$bjtu.edu.cn, 10010097@sxnu.edu.cn, orcawhale@163.com }
\begin{abstract}
In this paper, we construct an infinite families of group $G$  of order $4p^m$ which can be an automorphism group of some regular polytope with type $\{p, 2p\}$, where $m \geq 3$ and $p$ is an odd prime. For $p=3$, our polytopes are the regular toroidal maps $\{3, 6\}$.

%
\medskip
\f {\bf Keywords:} Regular $3$-polytopes, $p$-groups, automorphism groups.

\end{abstract}

\section{Introduction}

An \emph{abstract polytope} $\P$ is a partially-ordered set endowed with a rank function, satisfying certain conditions that arise naturally from a geometric setting. A
\emph{chain} of $\P$ is a totally ordered subset of $\P$, and a maximal chain is called a \emph{flag} of $\P$. An \emph{automorphism} of a polytope $\P$ is an incidence-preserving bijection from $\P$ to itself. The group of automorphisms of $\P$ is denoted by $\Gamma(\P)$.
The automorphism group of a polytope always acts freely on flags, and in the case where this action is transitive as well, the polytope is said to be \emph{regular}. The study of abstract regular polytopes has a rich history and been described comprehensively in the book by McMullen and Schulte~\cite{ARP}.

It is a natural question to try to describe all pairs $(\P, G)$, where $\P$ is a regular polytope and $G$ is the automorphism group
acting transitively on the flags of $\P$.
The atlas~\cite{atles1} contains information about all regular polytopes with automorphism group of order at most 2000.
 An interesting case is constituted by the pairs $(\P, G)$ with $G$ simple or almost simple.
The atlas~\cite{atles2} contains all regular polytopes whose automorphism group $G$ is an almost simple group such that $S \leq G\leq \Aut(S)$ and $S$ is a simple group of order less than one million,
and more striking results have been obtained for the symmetric groups $S_n$ and alternating groups $A_n$. Fernandes~{\em et al.}~\cite{fl,flm,sympolcorr} classified abstract regular polytopes of ranks $n-1$, $n-2$, $n-3$ and $n-4$ for $S_n$ respectively, and
Cameron {\em et al.}~\cite{CFLM2017} showed that the highest rank of an abstract regular polytope for $A_n$ is $\lfloor (n-1)/2 \rfloor$ when $n\geq 12$, which was known to be sharp by Fernandes~{\em et al.}~\cite{flm1,flm2}.
For a prime power $q$, more results about regular polytopes were obtained for the linear groups $\PSL(2, q)$ in \cite{psl2q}, $\PGL(2,q)$ in \cite{pgl2q}, $\PSL(3, q)$ and $\PGL(3, q)$ in \cite{psl3q}, $\PSL(4, q)$ in \cite{psl4q}, and for the Suzuki simple groups $\Sz(q)$ in \cite{suzuki}. Furthermore, Connor {\em et al.}~\cite{soclepsl2q} classified abstract regular polytopes for almost simple groups $G$ with $\PSL(2, q) \leq G \leq \PGaL(2, q)$. For more knowledge about this area, reader can refer the recent survey~\cite{Lee2019}.

Another interesting case is constituted by the pairs $(\P, G)$ with $G$ solvable, in this case, $G$ must have even order.
When it comes to solvable groups, except for abelian groups,
the first thing that come to mind are groups of orders $2^n$ or $2^np^m$, where $p$ is an odd prime. Schulte and Weiss~\cite{Problem} talked that the order $2^n$ and $2^np^m$ proved to be more difficult than others.

For $2$-groups, Conder~\cite{SmallestPolytopes} constructed an infinite families of regular polytopes of type $\{4, 4, \cdots, 4\}$ of order $2^n$, where $n$ depends on the rank.
The first author in this paper, Feng and Leemans \cite{HFL, HFL1} showed that all possible Schl\"afli types can be achieved for regular polytopes of order $2^n$ for $n \geq 5$. More recently,  Gomi {\em et al} \cite{sc1024} determined the non-degenerate regular polytopes of order $1024$.
For the groups of order $2^np$, some of small order and families of regular 3-polytopes arising from
the so-called tight regular $3$-polytopes, namely those with Schl\"afli type $\{k_1, k_2\}$ and order $2k_1k_2$, were obtained by Cunningham and Pellicer~\cite{GD2016}.
The first author in this paper, Feng and Leemans \cite{HFL2} constructed an infinite families of regular $3$-polytopes of order $2^np$ with type $\{2^{l_{1}}, 2^{l_{2}}p\}$ where $l_1, l_2 \geq 2$ and $l_1+l_2 \leq n-1$. They also given a family of regular $3$-polytopes of order $3\cdot2^n$ with type $\{6, 6\}$ for $n \geq 5$.

Naturally, we wish to construct regular polytopes $\P$ of order $2^np^m$ for $m \geq 2$.
Note that if $\P$ has rank at least $3$, then $4 \mid |\Aut(\P)|$.
Inspired by this, we
naturally assume that $|\Aut(\P)|=4p^m$ (now $\P$ has rank $3$). This paper
is devoted to provide a family of regular $3$-polytopes of order $4p^m$ with type $\{p, 2p\}$.

To state the result, we need to know some information about groups of order $p^m$. Let $G$ be a group and $H_p(G)$ the
subgroup of $G$ generated by all the elements of order different from $p$. $H_p(G)$ is called a {\em Hughes
subgroup} of $G$ (see~\cite{hughes,hughes1,hughes2}).
A group $G$ of order $p^m$ and nilpotency class $m-1$ is said to be maximal class if $m \geq 3$. The basic
material about these groups can be found in Blackburn~\cite{Blackburn} or Huppert~\cite[Chapter 3]{Huppert}.
If $G$ has an abelian maximal subgroup $A$, and $A=H_{p}(G)$, then, by \cite[Chapter 8, Example 8.3.3]{XQ}), 

$$G=\lg s_1, s_2, \cdots, s_r, s_{r+1}, \cdots, s_{p-1}, \b \rg\ \mbox{with}\ \mbox{the}\ \mbox{following}\ \mbox{relations}$$
$$ s_1^{p^e}, s_2^{p^e}, \cdots, s_r^{p^e}, s_{r+1}^{p^{e-1}}, \cdots s_{p-1}^{p^{e-1}}, \b^p,$$
$$[s_i, s_j],  \ {\rm for} \ 1 \leq i, j \leq p-1,\ \mbox{and}\ s_{k+1}=[s_k, \b],  \ {\rm for} \ 1 \leq k \leq p-2,$$
$$[s_{p-1}, \b]=s_1^{-\tbinom{p}{1}}s_2^{-\tbinom{p}{2}}\cdots s_{p-2}^{-\tbinom{p}{p-2}}s_{p-1}^{-\tbinom{p}{p-1}},$$
where $\tbinom{l_1}{l_2}$ denotes  the binomial coefficient indexed by $l_1$ and $l_2$, and
$$1 \leq e,1 \leq r \leq p-1,\ m=er+(e-1)(p-r-1)+1 \geq 3.$$
Moreover, for any $x \notin A=\lg s_1, s_2, \cdots, s_{p-1}\rg$, $x$ has order $p$.

Obviously, $G=\lg s_1, \b\rg$. Define
\begin{center}
\begin{tabular}{lll}
$\sigma:$ \qquad & $s_1 \mapsto s_1$ \quad &  and \quad $\beta \mapsto \beta^{-1}$,  \\
$\tau:$ \qquad & $s_1 \mapsto s_1^{-1}$ \quad & and \quad $\beta \mapsto \beta$.  \\
\end{tabular}
\end{center}
In this paper we prove the following theorem.

\begin{theorem}\label{maintheorem}
Let $R=G\rtimes \lg \sigma, \tau\rg$. Then $R$ can be the automorphism groups of a regular $3$-polytope of order $4p^m$ with type $\{p, 2p\}$, where $m \geq 3$ and $p$ is an odd prime.
\end{theorem}



\section{Background results}\label{backgroud}

\subsection{Abstract regular polytopes and String C-groups}

We know that abstract regular polytopes and string C-groups are the same mathematical objects. The link between them may be found for instance in~\cite[Chapter 2]{ARP}.
We take here the viewpoint of string C-groups because it is the easiest and most efficient one to define abstract regular polytopes.

Let $G$ be a group and $S=\{\rho_0,\cdots,\rho_{d-1}\}$ a generating set of involutions of $G$.
For $I\subseteq\{0,\cdots,d-1\}$, let $G_I$ denote the group generated by $\{\rho_i:i\in I\}$.
Suppose that
\begin{itemize}
\item[*]  for any $i,j\in \{0, \ldots, d-1\}$ with $|i-j|>1$, $\rho_i$ and $\rho_j$ commute (the \emph{string
property});
\item[*] for any $I,J\subseteq\{0,\cdots,d-1\}$,
$G_I\cap G_J=G_{I\cap J}\ \  (\mbox{the \emph{intersection property}})$.
\end{itemize}
Then the pair $(G,S)$ is called a {\em string C-group of rank} $d$ and the {\em order} of $(G,S)$ is simply the order of $G$. The {\em (Schl\"afli) type} of $(G,S)$ is the ordered set $\{p_1,\cdots, p_{d-1}\}$, where $p_i$ is the order of $\r_{i-1}\r_i$. If $(G,S)$  only satisfies the string property, then it is called a {\em string group generated by involutions} or \emph {sggi}. By the intersection property, $S$ is a minimal generating set of $G$.  It is known that string $C$-groups are same as automorphism groups of regular polytopes~\cite[Section 2E]{ARP}. In this paper we always assume that each $p_i$ is at least 3, otherwise, the generated group is a direct product of two smaller groups. If that happens, the string C-group (and the corresponding abstract regular polytope) is called {\em degenerate}.
The following proposition is
straightforward, and for details, one may see~\cite{MC}.
\begin{prop}\label{intersection}
The intersection property for a string $C$-group $(G,S)$ of rank $3$ is equivalent to that $S$ is a minimal generating set of $G$ and $\lg \r_0, \r_1\rg \cap \lg \r_1, \r_2\rg = \lg \r_1\rg $.
\end{prop}

\subsection{Group theory}

We use  standard notation for group theory, as in~\cite{GroupBook} for example.
In this subsection we briefly describe some knowledge of group theory we need.

\smallskip
Let $G$ be a group.   We define the {\em commutator\/} $[x, y]$ of elements $x$ and $y$ of $G$
by $[x, y]=x^{-1}y^{-1}xy$. The following result are elementary and so we give them without proof.

\begin{prop}\label{commutator}
Let $G$ be a group. Then, for any $x, y, z \in G$,
\begin{itemize}
\item [{\rm(1)}] $[xy, z]=[x, z]^y[y, z]$, $[x, yz]=[x, z][x, y]^z$;
\item [{\rm(2)}] $[x, y^{-1}]^y=[x, y]^{-1}$, $[x^{-1}, y]^x=[x, y]^{-1}$;
\item [{\rm(3)}] $[x^{-1}, y^{-1}]^{xy}=[x, y]$.
\end{itemize}
\end{prop}

\subsection{Combinatorial Enumeration}
In this subsection, we list some known results about combinatorial enumeration which be used in the sequel (see~\cite[Chapter 1]{MA}).

\begin{prop}\label{enumeration}
Let $n, m, k \in \mathbb{N}$. Then
\begin{itemize}
\item [\rm(1)] $\tbinom{n}{m} =\tbinom{n}{n-m}$;
\item [\rm(2)] $\tbinom{n+2}{m+1} =\tbinom{n+1}{m+1}+\tbinom{n+1}{m}$;
\item [\rm(3)] $\tbinom{n}{m}\tbinom{m}{k}=\tbinom{n}{k}\tbinom{n-k}{m-k}$;
\item [\rm(4)] $\sum\limits_{k=0}^m(-1)^{k}\tbinom{n}{k}=(-1)^m\tbinom{n-1}{m}$;
\end{itemize}
Furthermore, if $m <0$, then $\tbinom{n}{m}=0$.
\end{prop}

\section{The proof of Theorem 1.1 }

Assume $G$, $\sigma$ and $\tau$ as assumption above.

Let $s_{k+1}=[s_{k}, \b]$ for $k = p-1$ and $p$. Then
$s_p=[s_{p-1}, \b]=\prod\limits_{i=1}^{p-1}s_i^{-\tbinom{p}{i}}$, and
$$s_{p+1}=[s_p,\b]=s_p^{-1}s_p^{\b}=\prod\limits_{i=1}^{p-1}(s_i^{-1})^{-\tbinom{p}{i}}\prod\limits_{i=1}^{p-1}(s_i^{\beta})^{-\tbinom{p}{i}}$$
$$=\prod\limits_{i=1}^{p-1}(s_i^{-1}s_i^{\beta})^{-\tbinom{p}{i}}=\prod\limits_{i=1}^{p-1}[s_i, \beta]^{-\tbinom{p}{i}}=\prod\limits_{i=1}^{p-1}s_{i+1}^{-\tbinom{p}{i}}.$$

In order to prove Theorem 1.1, we need to prove the following two lemmas.

\begin{lem}\label{lemma1}
For any $1 \leq k \leq p$, we have
\begin{itemize}
\item [{\rm(1)}]$s_k^{\b^{1-k}}=\prod\limits_{i=k}^{p+1}s_i^{\tbinom{p+1-k}{i-k}}$;
\item [{\rm(2)}]$s_{k}^{\sigma}=\left\{
\begin{array}{ll}
s_k^{\b^{1-k}},   & k  \mbox{ is odd,} \\
(s_k^{\b^{1-k}})^{-1},  & k   \mbox{ is even.}
\end{array}
\right.$

\end{itemize}
\end{lem}

\demo
(1) Since $s_{k+1}=[s_k, \b]$, we have $s_k^{\b}=s_ks_{k+1}$.
Since $o(\b)=p$ and $k-1 \leq p$, we have $p-k+1 \geq 0$.
It follows that
$$s_k^{\b^{1-k}}=s_k^{\b^{p+1-k}}=(s_ks_{k+1})^{\b^{p-k}}=(s_ks_{k+1}^2s_{k+2})^{\b^{p-k-1}}=\cdots $$
$$=s_k^{\tbinom{p+1-k}{0}}s_{k+1}^{\tbinom{p+1-k}{1}}s_{k+2}^{\tbinom{p+1-k}{2}}\cdots s_{p}^{\tbinom{p+1-k}{p-k}}s_{p+1}^{\tbinom{p+1-k}{p+1-k}}.$$

\f (2) By induction on $k$. For $k=1$, we have $s_1^{\sigma}=s_1$. For $k=2$, by Proposition~\ref{commutator}, we have $$s_2^{\sigma}
=[s_1, \b]^{\sigma}=[s_1, \b^{-1}]=([s_1, \b]^{-1})^{\b^{-1}}=(s_2^{-1})^{\b^{1-2}}.$$

 Now this lemma is true for $k=1$ and $2$. By induction and Proposition~\ref{commutator},
if $k$ is odd, then $$s_{k}^{\sigma}=[s_{k-1}^{\sigma}, \b^{\sigma}]=[(s_{k-1}^{\b^{2-k}})^{-1}, \b^{-1}]=[s_{k-1}^{-1}, \b^{-1}]^{\b^{2-k}}=([s_{k-1}, \b]^{\b^{-1}s_{k-1}^{-1}})^{\b^{2-k}}.$$

Notice that $A=\lg s_1, s_2, \cdots, s_{p-1}\rg$ is an abelian maximal subgroup of $\mathcal{G}$ and hence $A \unlhd \mathcal{G}$.
It follows that $s_k \in A$ and $s_k^{\b} \in A$ for any $k \geq 1$. Thus,
$$s_{k}^{\sigma}=([s_{k-1}, \b]^{\b^{-1}s_{k-1}^{-1}})^{\b^{2-k}}=s_k^{(s_{k-1}^{-1})^{\b}\b^{1-k}}=s_k^{\b^{1-k}}.$$
If $k$ is even, then
$$s_{k}^{\sigma}=[s_{k-1}^{\sigma}, \b^{\sigma}]=[s_{k-1}^{\b^{2-k}}, \b^{-1}]=[s_{k-1}, \b^{-1}]^{\b^{2-k}}=(([s_{k-1}, \b]^{-1})^{\b^{-1}})^{\b^{2-k}}=(s_k^{\b^{1-k}})^{-1}.$$
\hfill\qed
\begin{lem}
$\sigma, \tau \in \Aut(G)$, and $\lg \sigma, \tau \rg \cong \mathbb{Z}_2 \times \mathbb{Z}_2$.
\end{lem}

\demo Obviously, $\sigma^2=\tau^2=1$ and $\sigma\tau=\tau\sigma$. Then $\lg \sigma, \tau\rg \cong \mathbb{Z}_2 \times \mathbb{Z}_2$. In the following we prove
$\sigma, \tau \in \Aut(G)$.

(1) We prove that $\sigma \in \Aut(G)$. Obviously, $\lg s_1^{\s}, \b^{\s}\rg=\lg s_1, \b^{-1}\rg=G$. Thus, we only need to show that the generating set $\{s_1^{\s}, s_2^{\s}, \cdots, s_{p-1}^{\s}, \b^{\s}\}$
of $G$ satisfy
the same relations as $\{s_1, s_2, \cdots, s_{p-1}, \b\}$.

Recall that $A=\lg s_1, s_2, \cdots, s_{p-1}\rg$ is abelian and $A \unlhd G$. By lemma~\ref{lemma1}(2),
$$s_{i}^{\sigma}=\left\{
\begin{array}{ll}
s_i^{\b^{1-i}},   & i \mbox{ is odd,} \\
(s_i^{\b^{1-i}})^{-1},  & i \mbox{ is even.}
\end{array}
\right.$$
Then $o(s_i)=o(s_i^{\s})$, $[s_i^{\s}, s_j^{\s}]=1$ for $1 \leq i, j \leq p-1$.
Also, it is easy to see that $s_{k+1}^{\sigma}=[s_k^{\s}, \b^{\s}]$ for $1\leq k \leq p-2$. Then we only need to show that
$$[s_{p-1}^{\sigma}, \b^{\sigma}]=(s_1^{\sigma})^{-\tbinom{p}{1}}(s_2^{\sigma})^{-\tbinom{p}{2}}(s_3^{\sigma})^{-\tbinom{p}{3}}\cdots (s_{p-1}^{\sigma})^{-\tbinom{p}{p-1}},$$
that is, $s_p^{\sigma}=\prod\limits_{i=1}^{p-1}(s_i^{\sigma})^{-\tbinom{p}{i}}.$

By Lemma~\ref{lemma1}(1), we have
$$s_1^{\sigma}=s_1, s_2^{\sigma}=\prod\limits_{i=2}^{p+1}s_i^{-\tbinom{p-1}{i-2}},
s_3^{\sigma}=\prod\limits_{i=3}^{p+1}s_i^{\tbinom{p-2}{i-3}}, \cdots, s_{p-1}^{\sigma}=
\prod\limits_{i=p-1}^{p+1}s_i^{-\tbinom{2}{i-(p-1)}}.$$

\f Then $\prod\limits_{i=1}^{p-1}(s_i^{\sigma})^{-\tbinom{p}{i}}
=s_1^{u_1}s_2^{u_2}\cdots s_{p+1}^{u_{p+1}}$, where $u_i$ comes from the following table

\setlength{\arraycolsep}{5pt}
\begingroup
\renewcommand*{\arraystretch}{1.8}
$$
 \begin{bmatrix}

   \begin{array}{c | c|cccccc}
*            &(s_1^{\sigma})^{-\tbinom{p}{1}}   &(s_2^{\sigma})^{-\tbinom{p}{2}}   &(s_3^{\sigma})^{-\tbinom{p}{3}}   &(s_4^{\sigma})^{-\tbinom{p}{4}}       &\cdots &\cdots &(s_{p-1}^{\sigma})^{-\tbinom{p}{p-1}}\\ \hline
u_1       &-\tbinom{p}{1}      & 0                                                & 0                                                     &0 &\cdots &\cdots         &0                               \\ \hline
u_2       &0&\tbinom{p-1}{0}\tbinom{p}{2}      & 0                                                & 0                                                     &\cdots &\cdots &0                                      \\
u_3        &0&\tbinom{p-1}{1}\tbinom{p}{2}      & -\tbinom{p-2}{0}\tbinom{p}{3}      & 0                                                    &\cdots &\cdots &0                                       \\
u_4        &0&\tbinom{p-1}{2}\tbinom{p}{2}      & -\tbinom{p-2}{1}\tbinom{p}{3}      & \tbinom{p-3}{0}\tbinom{p}{4}        &\cdots &\cdots &0                                        \\
u_5      &0&\tbinom{p-1}{3}\tbinom{p}{2}      & -\tbinom{p-2}{2}\tbinom{p}{3}      & \tbinom{p-3}{1}\tbinom{p}{4}        &\cdots &\cdots &0                                     \\
\vdots    &\vdots&\vdots                                          & \vdots                                         & \vdots                                            &\cdots &\cdots &  \vdots                             \\
\vdots    &\vdots&\vdots                                          & \vdots                                         & \vdots                                            &\cdots &\cdots &   \vdots                                \\
u_{p-1}  &0&\tbinom{p-1}{p-3}\tbinom{p}{2}   & -\tbinom{p-2}{p-4}\tbinom{p}{3}   & \tbinom{p-3}{p-5}\tbinom{p}{4}     &\cdots &\cdots &\tbinom{2}{0}\tbinom{p}{p-1}   \\ \hline
u_p     &0&\tbinom{p-1}{p-2}\tbinom{p}{2}   & -\tbinom{p-2}{p-3}\tbinom{p}{3}   & \tbinom{p-3}{p-4}\tbinom{p}{4}     &\cdots &\cdots &\tbinom{2}{1}\tbinom{p}{p-1} \\ \hline
u_{p+1}&0&\tbinom{p-1}{p-1}\tbinom{p}{2}   & -\tbinom{p-2}{p-2}\tbinom{p}{3}   & \tbinom{p-3}{p-3}\tbinom{p}{4}     &\cdots &\cdots &\tbinom{2}{2}\tbinom{p}{p-1}

  \end{array}

  \end{bmatrix}
$$
\endgroup

Obviously, $u_1=-\tbinom{p}{1}$. By Proposition~\ref{enumeration}(1), we have
$$u_{p+1}=\sum\limits_{j=2}^{p-1}(-1)^{j}\tbinom{p}{j}=\tbinom{p}{p-1}+
\sum\limits_{j=2}^{\frac{p-3}{2}}\left[(-1)^{j}\tbinom{p}{j}+(-1)^{p-j}\tbinom{p}{p-j}\right]=\tbinom{p}{p-1}=p.$$
\f For $2 \leq i \leq p-1$, we have $u_i=\sum\limits_{j=2}^{i}(-1)^j\tbinom{p+1-j}{i-j}\tbinom{p}{j}$.
For $i=p$, we have $$u_p-\tbinom{1}{0}\tbinom{p}{p}=u_p-1=
\sum\limits_{j=2}^{p}(-1)^j\tbinom{p+1-j}{p-j}\tbinom{p}{j}.$$

By Proposition~\ref{enumeration}$(2)$$(3)$, we have
$$\tbinom{p+1-j}{i-j}\tbinom{p}{j}=\left[\tbinom{p-j}{i-j}+\tbinom{p-j}{i-j-1}\right]\tbinom{p}{j}
=\tbinom{p-j}{i-j}\tbinom{p}{j}+\tbinom{p-j}{i-j-1}\tbinom{p}{j}=
\tbinom{p}{i}\tbinom{i}{j}+\tbinom{p}{i-1}\tbinom{i-1}{j}.$$
Together with Proposition~\ref{enumeration}(4), for $2 \leq i \leq p$, we have
$$\sum\limits_{j=2}^{i}(-1)^j\left[\tbinom{p}{i}\tbinom{i}{j}+\tbinom{p}{i-1}\tbinom{i-1}{j}\right]$$
$$=\tbinom{p}{i}\left[i-1+\sum\limits_{j=0}^{i}(-1)^j\tbinom{i}{j}\right]+
\tbinom{p}{i-1}\left[i-2+\sum\limits_{j=0}^{i}(-1)^j\tbinom{i-1}{j}\right]$$
$$=\tbinom{p}{i}\left[i-1+(-1)^i\tbinom{i-1}{i}\right]+
\tbinom{p}{i-1}\left[i-2+(-1)^i\tbinom{i-2}{i}\right]$$
$$=\tbinom{p}{i}(i-1)+\tbinom{p}{i-1}(i-2)$$
$$=-\tbinom{p}{i}+(p-i+1)\tbinom{p}{i-1}+(i-2)\tbinom{p}{i-1}$$
$$=-\tbinom{p}{i}
+(p-1)\tbinom{p}{i-1},$$
and hence
$$u_p=-\tbinom{p}{p}
+(p-1)\tbinom{p}{p-1}+1=(p-1)\tbinom{p}{p-1},$$
$$u_i=-\tbinom{p}{i}
+(p-1)\tbinom{p}{i-1}, \quad 2 \leq i \leq p-1.$$

By Lemma~\ref{lemma1} and $s_{p+1}=[s_p,\b]$, we have $s_p^{\sigma}=s_p^{\b^{1-p}}=s_p^{\b}=s_ps_{p+1}.$
The result is now immediate since $$s_1^{u_1}s_2^{u_2}\cdots s_p^{u_p}s_{p+1}^{u_{p+1}}$$
$$=s_1^{-\tbinom{p}{1}}\left[\prod\limits_{i=2}^{p-1}s_i^{-\tbinom{p}{i}
+(p-1)\tbinom{p}{i-1}}\right]s_{p}^{(p-1)\tbinom{p}{p-1}}s_{p+1}^{p}$$
$$=\left[\prod\limits_{i=1}^{p-1}s_i^{-\tbinom{p}{i}}\right]\left[(\prod\limits_{i=2}^{p}s_i^{-\tbinom{p}{i-1}})^{1-p}\right]s_{p+1}^p$$
$$=s_ps_{p+1}^{1-p}s_{p+1}^{p}=s_{p}s_{p+1}=s_p^{\sigma}.$$

(2) We prove that $\tau \in \Aut(G)$. Since $A$ is abelian and $s_{i+1}=[s_i, \b]$ for $1 \leq i \leq p-2$, by Proposition~\ref{commutator}, we have
$$s_{i+1}^{\tau}=[s_i, \b]^{\tau}=[s_i^{-1}, \b]=(([s_i, \b])^{-1})^{s_i^{-1}}=(s_{i+1}^{-1})^{s_i^{-1}}=s_{i+1}^{-1}.$$

We need to show that the generating set $\{s_1^{-1}, s_2^{-1}, \cdots, s_{p-1}^{-1}, \b\}$ of $G$ satisfy
the same relations as $\{s_1, s_2, \cdots, s_{p-1}, \b\}$.
Obviously,  $o(s_i)=o(s_i^{\tau})$, $[s_i^{\tau}, s_j^{\tau}]=1$ for $1 \leq i, j \leq p-1$.
Also, we have $$[s_k, \b]^{\tau}=[s_k^{\tau}, \b^{\tau}]=[s_{k}^{-1},\ \b]=([s_k, \b]^{-1})^{s_k^{-1}}=s_{k+1}^{-1}=s_{k+1}^{\tau}\ \mbox{for}\ 1\leq k \leq p-2.$$

Now we only need to show that
$$[s_{p-1}^{\tau}, \b^{\tau}]=(s_1^{\tau})^{-\tbinom{p}{1}}(s_2^{\tau})^{-\tbinom{p}{2}}(s_3^{\tau})^{-\tbinom{p}{3}}\cdots (s_{p-1}^{\tau})^{-\tbinom{p}{p-1}}.$$
It is easy to see that $$[s_{p-1}^{\tau}, \b^{\tau}]=[s_{p-1}^{-1}, \b]=([s_{p-1}, \b]^{-1})^{s_{p-1}^{-1}}=(s_{p}^{-1})^{s_{p-1}^{-1}}=s_{p}^{-1}.$$ On the other hand,
$$\prod\limits_{i=1}^{p-1}(s_i^{\tau})^{-\tbinom{p}{i}}=\prod\limits_{i=1}^{p-1}(s_i^{-1})^{-\tbinom{p}{i}}
=(\prod\limits_{i=1}^{p-1}s_i^{-\tbinom{p}{i}})^{-1}=s_p^{-1}.$$
It follows that $\tau \in \Aut(G)$.
\hfill\qed

Now we are ready to prove Theorem~\ref{maintheorem}.

\medskip
\f{\bf Proof of Theorem 1.1.}
Let $\r_0=s_1\tau\sigma, \r_1=\b\tau\sigma, \r_2=\sigma$. Since $o(\b)=p, o(\tau)=2$ and $[\b, \tau]=1$, we have
$(\r_1\r_2)^p=(\b\tau)^p=\tau$, and hence $\b=(\r_1\r_2)^{p+1}$. It follows that $R=\lg \r_0, \r_1, \r_2\rg$.

Note that
$$\r_0^2=s_1\tau\sigma s_1\tau\sigma=s_1s_1^{\tau\sigma}=s_1s_1^{-1}=1,$$
$$\r_1^2=\b\tau\sigma\b\tau\sigma=\b\b^{\sigma\tau}=\b\b^{-1}=1,$$
$$(\r_0\r_2)^2=(s_1\tau\sigma\sigma)^2=s_1s_1^{\tau}=s_1s_1^{-1}=1.$$
Then $\{R, \{\r_0, \r_1, \r_2\}\}$ is an {\em sggi}. We shall prove $\{R, \{\r_0, \r_1, \r_2\}\}$ is a string C-group.

Note that for any $x \notin \lg s_1, s_2, \cdots, s_{p-1}\rg$, we have $o(x)=p$.
Since $$\r_0\r_1=s_1\tau\sigma\b\tau\sigma=s_1\b^{-1} \notin \lg s_1, s_2, \cdots, s_{p-1}\rg,$$  $o(\r_0\r_1)=p$.
Since $(\r_1\r_2)^2=(\b\tau\sigma\tau)^2=(\b\sigma)^2=\b^2$, $o(\r_1\r_2)=2p$.
If follows that $|\lg \r_0, \r_1\rg|=2p$ and $|\lg \r_1, \r_2\rg|=4p$.

Finally, if $\lg \r_0, \r_1\rg \cap \lg \r_1, \r_2\rg >\lg \r_1\rg$, then $\lg \r_0, \r_1\rg \leq \lg \r_1, \r_2\rg$.
It means that $G=\lg \r_1, \r_2\rg$, which is impossible because $|G|=4p^m>4p$.
Thus, $\lg \r_0, \r_1\rg \cap \lg \r_1, \r_2\rg=\lg \r_1\rg$.
By Proposition~\ref{intersection}, $\{R, \{\r_0, \r_1, \r_2\}\}$ is a string C-group
with type $\{p, 2p\}$. That means  $R$ can be the automorphism groups of a regular $3$-polytope of order $4p^m$ with type $\{p, 2p\}$ for $m \geq 3$.
\hfill\qed

\f {\bf Acknowledgements:} This work was supported by the National Natural Science Foundation of China 11771258 and 11971280.

\end{document}